\documentclass[12pt]{article}
\usepackage{epic,eepic,epsf,epsfig}
\usepackage{amsmath}
\usepackage{amssymb}
\usepackage{amsbsy}

\usepackage{amsfonts}
\newtheorem{theorem}{Theorem}[section]%
\newtheorem{cor}[theorem]{Corollary}%
\newtheorem{remark}[theorem]{Remark}%

\def\f{\noindent}

\newcommand{\qed}{\mbox{\raisebox{0.7ex}{\fbox{}}} \vspace{4truemm}}
\def\demo{\f {\bf Proof.}\hskip10pt}

\begin{document}

\baselineskip 16pt

\title{ \vspace{-1.2cm}
On a finite group in which every non-supersolvable maximal subgroup is a TI-subgroup
\thanks{\scriptsize This work was supported in part by NSFC (12571024) and Shandong Provincial Natural Science Foundation, China (ZR2017MA022).
\newline
 \hspace*{0.5cm} \scriptsize $^{\ast\ast}$Corresponding
  author.\newline
 \hspace*{0.5cm} \scriptsize{E-mail addresses:} jtshi2005@ytu.edu.cn\,(J. Shi),\,lifh2026@s.ytu.edu.cn\,(F. Li).}}

\author{Jiangtao Shi\,$^{\ast\ast}$,\,\,Fanghua Li\\
\\
{\small{\em School of Mathematics and Information Sciences,
Yantai University, Yantai 264005, P.R. China}}}

\date{ }

\maketitle \vspace{-.8cm}

\begin{abstract}
We provide a complete characterization of a finite group in which every non-supersolvable maximal subgroup is a TI-subgroup, which generalizes some known related results.

\medskip
\f {\bf Keywords:} non-supersolvable maximal subgroup; TI-subgroup; Frobenius group; solvable group
\end{abstract}

\section{Introduction}\label{s1}

In this paper all groups are assumed to be finite. Let $G$ be a group and $H$ a subgroup of $G$. Then $H$ is said to be a TI-subgroup of $G$ if $H^g\cap H$=1 or $H$ for each $g\in G$,
that is $H^g\cap H=1$ for each $g\in G\backslash N_G(H)$. Investigating the structure of a group in which some particular maximal subgroups are TI-subgroups, Lu, Pang and Zhong {\rm\cite[Theorem 3.5]{lu2013}}
showed that every non-nilpotent maximal subgroup of a group $G$ is a TI-subgroup if and only if every non-nilpotent maximal subgroup of $G$ is normal. Note that Shi, Li and Shen
{\rm\cite[Theorem 1.7]{shilishen}} proved that if every non-nilpotent maximal subgroup of a group $G$ is normal then $G$ is either $q$-nilpotent
or $q$-closed for each prime divisor $q$ of $|G|$. Furthermore, Shi and Xu {\rm\cite{shixu2025}} provided a complete characterization of a group
in which every non-nilpotent maximal subgroup is normal.

As a generalization, Shi {\rm\cite{s6}} gave a complete description of a group $G$ in which every non-nilpotent maximal subgroup of order
divisible by $p$ is a TI-subgroup for any prime divisor $p$ of $|G|$. Combining the coprime action of groups, Lu, Li, Zhang et al. {\rm\cite{lu2023}}
characterized a group in which every non-abelian maximal invariant subgroup is a TI-subgroup or a subnormal subgroup. In {\rm\cite{xushishan}}
and {\rm\cite{shishan}} a group in which every non-nilpotent maximal invariant subgroup of order divisible by $p$ is a TI-subgroup
and a group in which every non-abelian maximal invariant subgroup of order divisible by $p$ is a TI-subgroup are investigated respectively, where $p$ is any prime divisor of the order of
the group.

As a further generalization of above research, considering the group in which every non-supersolvable maximal subgroup is a TI-subgroup, we have the following result,
the proof of which is given in Section~\ref{s2}.

\begin{theorem}\ \ \label{th1} Suppose that $G$ is a group. Then every non-supersolvable maximal subgroup of $G$ is a TI-subgroup if and only of
one of the following statements holds.

$(1)$ every non-supersolvable maximal subgroup of $G$ is normal;

$(2)$ $G={Z_p}^m\rtimes(Q_8\rtimes H)$ is a Frobenius group with kernel ${Z_p}^m$ and complement $Q_8\rtimes H$, where $m>1$, $p$ is the largest prime divisor of $|G|$, $H$ is a cyclic group of odd order, $Q_8\rtimes H$ acts irreducibly on ${Z_p}^m$ and $[Q_8,H]\neq 1$,
every non-cyclic maximal subgroup of  $Q_8\rtimes H$ contains $Q_8$, every cyclic maximal subgroup of $Q_8\rtimes H$ is a conjugate of $Z(Q_8)\times H$ and ${Z_p}^m\rtimes(Z(Q_8)\times H)$ is supersolvable.
\end{theorem}

The following corollary is a direct consequence of Theorem~\ref{th1}.

\begin{cor}\ \ \label{c1} Suppose that $G$ is a non-solvable group. Then every non-supersolvable maximal subgroup of $G$ is a TI-subgroup if and only if every non-supersolvable maximal subgroup of $G$ is normal.
\end{cor}

\begin{remark}\ \ \label{r1} {\rm In Theorem~\ref{th1} ${Z_p}^m\rtimes(Z(Q_8)\times H)$ is supersolvable implies that $Z(Q_8)\times H\lesssim Z_{p-1}$.}
\end{remark}

\begin{remark}\ \ \label{r2} {\rm Let $G\cong(Z_7\times Z_7)\rtimes(Q_8\rtimes Z_3)\cong (Z_7\times Z_7)\rtimes\rm SL(2,3)$, where $Q_8\rtimes Z_3$ is a subgroup of ${\rm Aut}(Z_7\times Z_7)\cong\rm GL(2,7)$
and $Q_8\rtimes Z_3$ acts irreducibly on $Z_7\times Z_7$.
It is easy to see that every non-supersolvable maximal subgroup of $G$ is either normal or a Frobenius complement of $G$ which is a conjugate of $Q_8\rtimes Z_3$.}
\end{remark}

\begin{remark}\ \ \label{r3} {\rm Note that a group in which every non-supersolvable maximal subgroup is normal might not be solvable. For example, let $G\cong{\rm Aut}(\rm SL(3,2))$. Then
${\rm Aut}(\rm SL(3,2))$ has a unique non-supersolvable maximal subgroup which is isomorphic to $\rm SL(3,2)$ and $\rm SL(3,2)$ is non-solvable.}
\end{remark}

\begin{remark}\ \ \label{r4} {\rm If every non-supersolvable maximal subgroup of a group $G$ is normal, we cannot get that $G$ is either $q$-nilpotent
or $q$-closed for each prime divisor $q$ of $|G|$. For example, let $G\cong S_4$, the symmetric group of degree 4. It is easy to see that $G$ is neither $q$-nilpotent nor $q$-closed whenever $q=2$ or 3.}
\end{remark}

\begin{remark}\ \ \label{r5} {\rm It is easy to show that a group in which every non-supersolvable subgroup is normal is solvable by analyzing the counterexample of the smallest order. Then
a group in which every non-supersolvable subgroup is a TI-subgroup is solvable combining Theorem~\ref{th1}.}
\end{remark}

\section{Proof of Theorem~\ref{th1}}\label{s2}

\demo {\bf Part I}: In this part, we prove the necessity part.

If $G$ has non-normal non-supersolvable maximal subgroups. Let $M$ be a non-normal non-supersolvable maximal subgroup of $G$. Then $M=N_G(M)$. Since $M$ is a TI-subgroup of $G$ by the hypothesis, one has $M\cap M^g=1$ for each $g\in G\backslash N_G(M)=G\backslash M$. It follows that $G$ is a Frobenius group with $M$ being a complement. Let $N$ be the kernel
of $G$. Then $G=N\rtimes M$. Since $M$ is maximal in $G$, $N$ is a minimal normal subgroup of $G$, which implies that $M$ acts irreducibly on $N$. Note that $N$ is nilpotent by {\rm\cite[Theorem 10.5.6]{rob}}.
Then $N$ is an elementary abelian subgroup of $G$ of prime-power order. Assume $N\cong{Z_p}^m$, where $p$ is a prime divisor of $|G|$ and $m\geq 1$.

By {\rm\cite[Theorem 10.5.6]{rob}}, every Sylow $p$-subgroup of $M$ is cyclic if $p>2$ and cyclic or a generalized quaternion group if $p=2$. Since a group in which every Sylow subgroup is cyclic is supersolvable
and $M$ is non-supersolvable, $M$ has even order and the Sylow 2-subgroup of $M$ is a generalized quaternion group. Let $R$ be a Sylow 2-subgroup of $M$. It is obvious that $M>R$.

Let $M_0$ be any maximal subgroup of $M$. Next we show that $N\rtimes M_0$ is supersolvable if and only if $M_0$ is cyclic.

First suppose that $M_0$ is cyclic. If $N\rtimes M_0$ is non-supersolvable. Then $N\rtimes M_0$ is a TI-subgroup of $G$. Since $(N\rtimes M_0)^g\cap(N\rtimes M_0)=(N^g\rtimes {M_0}^g)
\cap(N\rtimes M_0)=(N\rtimes {M_0}^g)\cap(N\rtimes M_0)\geq N\neq 1$ for each $g\in G$, $N\rtimes M_0$ is normal in $G$. Then $M_0=(N\cap M)\rtimes M_0=(N\rtimes M_0)\cap M$ is normal in $M$.
It follows that $M$ is supersolvable, a contradiction. Therefore, $N\rtimes M_0$ is supersolvable.

Next suppose that $N\rtimes M_0$ is supersolvable. It implies that there exists a minimal normal subgroup $N_0$ of $N\rtimes M_0$ of prime order such that $N_0\leq N$. Then $N_0\cong Z_p$.
Considering the group $N_0\rtimes M_0$, which is also a Frobenius group. By $\rm N/C$-theorem, one has $M_0\cong M_0/(N_0\cap M_0)\cong N_0\rtimes M_0/N_0=N_{N_0\rtimes M_0}(N_0)/C_{N_0\rtimes M_0}(N_0)\lesssim{\rm Aut}(N_0)
\cong Z_{p-1}$. Then $M_0$ is cyclic.

If every maximal subgroup of $M$ is cyclic. Then $M$ is either a cyclic group or a minimal non-cyclic group, which is supersolvable, a contradiction. If every maximal subgroup of $M$ is non-cyclic,
arguing as above, one has that every maximal subgroup of $M$ is normal, which implies that $M$ is nilpotent, also a contradiction. Therefore, $M$ not only has cyclic maximal subgroups but also has non-cyclic
maximal subgroups. It follows that $M$ is solvable and then $G$ is solvable.

Claim that $R$ is normal in $M$. Otherwise, assume that $R$ is not normal in $M$. Then $N_M(R)<M$. Let $M_1$ be a maximal subgroup of $M$ such that $N_M(R)\leq M_1$. Note that $M_1$ is non-cyclic.
Then $M_1$ is normal in $M$. By the Frattini argument, one has $M=M_1N_M(R)=M_1$, a contradiction. Thus $N_M(R)=M$, which implies that $R$ is normal in $M$. By Schur-Zassenhaus theorem, let $H$ be a complement
of $R$ in $M$. Then $M=R\rtimes H$, where $H>1$.

Let $H_0$ be any maximal subgroup of $H$. Then $N\rtimes(R\rtimes H_0)$ is a non-supersolvable maximal subgroup of $G$. Arguing as above, one has that $H_0$ is normal in $H$. Then $H$ is nilpotent. It follows that
$H$ is a cyclic group of odd order since every Sylow subgroup of $H$ is cyclic. Obviously $[R,H]\neq 1$.

Assume $R=Q_{4n}=\langle a,\,b|\,a^{2n}=b^4=1,\,b^{-1}ab=a^{-1}\rangle$. If $|R|\geq 16$. Then $\langle a\rangle$ is a characteristic subgroup of $R$. Then $\langle a\rangle$ is also normal in $M$.
Considering the subgroup $\langle a\rangle\rtimes H$. Since $\langle a\rangle$ is cyclic and 2 is the smallest prime divisor of $|\langle a\rangle\rtimes H|$. It follows that $\langle a\rangle\rtimes H$ is
2-nilpotent, which implies that $\langle a\rangle\rtimes H=\langle a\rangle\times H$ is cyclic. Furthermore, $\langle a\rangle\times H$ is a normal cyclic maximal subgroup of $R\rtimes H$ since $|R\rtimes H:\langle a\rangle\times H|=2$. It follows that $R\rtimes H$ is supersolvable, a contradiction. Therefore, $|R|=8$ and then $R=Q_8$.

Claim that every non-cyclic maximal subgroup of $M$ contains $Q_8$. Let $K$ be any non-cyclic maximal subgroup of $M$. If $Q_8\nleq K$. Then $Q_8\cap K<Q_8$ and $Q_8\cap K$ is normal in $M$
since both $Q_8$ and $K$ are normal in $M$. Let $K_0$ be a $2'$-Hall subgroup of $K$. Then $M=Q_8K=Q_8\rtimes K_0=Q_8\rtimes H$. It follows that $K_0\cong Q_8K_0/Q_8=Q_8H/Q_8\cong H$ is cyclic.
Note that $K=K\cap M=K\cap(Q_8\rtimes K_0)=(Q_8\cap K)\rtimes K_0$. Since $Q_8\cap K<K$, $Q_8\cap K$ is cyclic. It implies that $(Q_8\cap K)\rtimes K_0$ is 2-nilpotent since 2 is the smallest prime divisor of
$|(Q_8\cap K)\rtimes K_0|$. Thus $K=(Q_8\cap K)\rtimes K_0=(Q_8\cap K)\times K_0$ is cyclic, a contradiction. Therefore, every non-cyclic maximal subgroup of $M$ contains $Q_8$.

It is easy to see that $Z(Q_8)\rtimes H=Z(Q_8)\times H$ is cyclic. Moreover,  $Z(Q_8)\times H$ is a cyclic maximal subgroup of $M$ since $M$ is non-supersolvable which implies that the index of any cyclic maximal subgroup
in $M$ is greater than 2.

Let $T$ be any cyclic maximal subgroup of $Q_8\rtimes H$. Here $Q_8\cap T$ is a normal Sylow 2-subgroup of $T$. Let $T_0$ be a $2'$-Hall subgroup of $T$. Then $T=(Q_8\cap T)\rtimes T_0$. One has $M=Q_8T=Q_8\rtimes T_0$, where
both $T_0$ and $H$ are complements of $Q_8$ in $M$. Then $T_0$ and $H$ are conjugate in $M$. Since $T$ is a cyclic maximal subgroup of $M$ and $M$ is non-supersolvable, one has $2\leq |Q_8\cap T|<4$. Then $Q_8\cap T=Z(Q_8)$.
It follows that $T=Z(Q_8)\times T_0$ and $Z(Q_8)\times H$ are conjugate in $M$.

Note that $N\rtimes(Z(Q_8)\times H)$ is supersolvable. Let $N_1$ be a normal subgroup of $N\rtimes(Z(Q_8)\times H)$ of prime order satisfying $N_1\leq N$. Then $N_1\cong Z_p$. Considering the group
$Z_p\rtimes(Z(Q_8)\times H)$. By N/C-theorem, one has $Z(Q_8)\times H\cong Z_p\rtimes(Z(Q_8)\times H)/Z_p=N_G(Z_p)/C_G(Z_p)\lesssim{\rm Aut}(Z_p)=Z_{p-1}$. It follows that $p$ is the largest prime
divisor of $|G|$.

{\bf Part II}: In this part, we prove the sufficiency part.

If $G$ is a group belonging to case (1). It is obvious that every non-supersolvable maximal subgroup of $G$ is a TI-subgroup. In the following assume that $G$ is a group belonging to case (2).
Let $L$ be any non-supersolvable maximal subgroup of $G$.

Case $(i)$: Assume $N\leq L$. Then $L=L\cap G=L\cap(N\rtimes(Q_8\rtimes H))=N\rtimes(L\cap(Q_8\rtimes H))$.
Let $L_0=L\cap(Q_8\rtimes H)$. It is
obvious that $L_0$ is a maximal subgroup of $Q_8\rtimes H$. If $Q_8\nleq L_0$. By the hypothesis, $L_0$ is cyclic, $L_0$ is a conjugate of $Z(Q_8)\times H$ and $N\rtimes L_0$ is supersolvable,
that is $L$ is supersolvable, a contradiction. Therefore, $Q_8\leq L_0$. Then $L=N\rtimes(Q_8\rtimes(L\cap H))$. Since $H$ is cyclic, $L\cap H$ is normal in $H$. It follows that $L$ is normal in $G$,
which implies that $L$ is a TI-subgroup of $G$.

Case $(ii)$: Assume $N\nleq L$. Then $G=NL$. Since $N\cap L$ is not only normal in $L$ but also normal in $N$, $N\cap L$ is normal in $G$. It is obvious that $N\cap L<N$. As $N$ being a minimal normal
subgroup of $G$ by the hypothesis, one has $N\cap L=1$. Then $G=N\rtimes L$. Note that both $Q_8\rtimes H$ and $L$ are complements of $N$ in $G$. Moreover, one has that $G$ is solvable by the hypothesis.
It follows that $Q_8\rtimes H$ and $L$ are conjugate in $G$ by Schur-Zassenhaus theorem. Then $L$ is also a Frobenius complement of $G$, which implies that $L$ is a TI-subgroup of $G$.
\hfill\qed

\end{document}